\documentclass[12pt,a4paper]{amsart}

\usepackage{amssymb}

\usepackage{Alegreya}
\usepackage[libertine]{newtxmath}
\usepackage[headings]{fullpage}
\usepackage{paralist}

\usepackage{color}
\usepackage[normalem]{ulem}

\usepackage[colorlinks=true,linkcolor=blue, citecolor=blue, urlcolor=blue,%
pagebackref=true]{hyperref}

\makeatletter

\@addtoreset{equation}{section}
\def\theequation{\thesection.\@arabic \c@equation}

\def\@citecolor{blue}
\def\@linkcolor{blue}
\def\@urlcolor{blue}
\def\theenumi{\@alph\c@enumi}

\makeatother

\theoremstyle{plain}
\newtheorem{theorem}[equation]{Theorem}
\newtheorem{lemma}[equation]{Lemma}
\newtheorem{corollary}[equation]{Corollary}
\newtheorem{proposition}[equation]{Proposition}

\theoremstyle{definition}

\newtheorem{remark}[equation]{Remark}

\newtheorem{example}[equation]{Example}

\newtheorem{definition}[equation]{Definition}

\newtheorem{notation}[equation]{Notation}

\newtheorem{discussion}[equation]{Discussion}
\newenvironment{discussionbox}[1][]{%
    \begin{discussion}[#1]\pushQED{\qed}}{\popQED \end{discussion}}

\newtheorem{observation}[equation]{Observation}

\newtheorem{construction}[equation]{Construction}

\newcommand{\calB}{\mathcal B}

\newcommand{\bbF}{\mathbb F}

\DeclareMathOperator{\Sym}{Sym}
\newcommand{\naturals}{\mathbb{N}}
\newcommand{\ints}{\mathbb{Z}}

\def\to{\longrightarrow}
\DeclareMathOperator{\Trace}{Tr}

\DeclareMathOperator{\height}{ht}
\DeclareMathOperator{\codim}{codim}
\DeclareMathOperator{\charact}{char}

\newcommand{\define}[1]{\emph{#1}}
\DeclareMathOperator{\image}{Im}
\DeclareMathOperator{\GL}{GL}

\makeatletter
\def\RDerChar{\mathbf{R}}
\def\RDer{\@ifnextchar[{\R@Der}{\ensuremath{\RDerChar}}}
\def\R@Der[#1]{\ensuremath{\RDerChar^{#1}}}

\newif\ifreadkumminibib
\readkumminibibfalse

\makeatother

\newcommand{\hilbertIdeal}{\mathfrak{h}}

\title{On Hilbert ideals for a class of $p$-groups in characteristic $p$}
\author{Manoj Kummini}
\address{Chennai Mathematical Institute, Siruseri, Tamilnadu 603103. India}
\email{mkummini@cmi.ac.in}

\author{Mandira Mondal}
\address{Chennai Mathematical Institute, Siruseri, Tamilnadu 603103. India}
\email{mandiram@cmi.ac.in}

\thanks{Authors were supported by the grant CRG/2018/001592 from Science
and Engineering Research Board, India and by an Infosys Foundation
fellowship.}
\thanks{MM thanks National Board of Higher Mathematics.}

\dedicatory{In memory of Prof. C.~S.~Seshadri}
\begin{document}

\begin{abstract}
Let $p$ be a prime number, $\Bbbk$ a field of characteristic $p$
and $G$ a finite $p$-group.
Let $V$ be a finite-dimensional linear representation of $G$ over $\Bbbk$.
Write $S = \Sym V^*$.
For a class of $p$-groups which we call generalised Nakajima groups, we
prove the following:
\begin{enumerate}
\item The Hilbert ideal is a complete intersection.
As a consequence, for the case of generalised Nakajima groups, we prove a
conjecture of Shank and Wehlau (reformulated by Broer) that asserts that
if the invariant subring $S^G$ is a direct summand of $S$ as $S^G$-modules
then $S^G$ is a polynomial ring.
\item
The Hilbert ideal has a generating set with elements of degree at most $|G
|$.
This bound is conjectured by Derksen and Kemper.
\end{enumerate}
\end{abstract}

\maketitle

\section{Introduction}

Let $\Bbbk$ be a field and $V$ a finite dimensional $\Bbbk$-vector-space.
Let $G\subseteq \GL(V)$ be finite group.
Then the action of $G$ on $V$ induces an action on $V^*$.
Namely, if $\{v_1,\ldots,v_n\}$ is a basis of $V$ with
$\{x_1,\ldots,x_n\}\subseteq V^*$ its dual basis, then for each $g\in G$,
one sets $g x_i(v)=x_i(g^{-1}v)$ for all $v\in V$ and for all
$i=1,\ldots,n$.
This action extends to a graded $\Bbbk$-algebra automorphism of
$S :=\Sym V^*=\Bbbk[x_1,\ldots,x_n]$.
The \define{Hilbert ideal}
$\hilbertIdeal_{G, S}$ is $(S^G)_+S$, i.e, the $S$-ideal 
generated by all invariants of positive degree.
This is an important object in invariant theory; in some cases, the
generators of $\hilbertIdeal_{G, S}$ as an $S$-ideal generate $S^G$ as a
$\Bbbk$-algebra; see,
e.g.,~\cite[2.2.10]{DerksenKemperComputationalInvariantThy2ed2015}.

Throughout this article, $\charact \Bbbk  = p > 0$ and $G$ is a $p$-group.
By \cite[Proposition~4.0.2]{CampbellWehlauModularInvThy11}, there exists a
basis $\mathcal{B}^*=\{v_1,\ldots,v_n\}$ of $V$ such that for all $g\in G$,
and for all $i=1,\ldots,n$, $gv_i-v_i\in \Bbbk\langle v_{i+1},\ldots,
v_n\rangle$.
Moreover, if $\mathcal{B}=\{x_1,\ldots,x_n\}$ is the dual basis of
$\mathcal{B}^*$, then for all $g\in G$, and for all $i=1,\ldots,n$,
$$
gx_i-x_i\in \Bbbk\langle x_{1},\ldots, x_{i-1}\rangle.
$$
For $g\in G$, define $\beta_g=\text{max}\{j\mid x_j \text{ appears in }
gx_i-x_i \text{ for some }i\}$.

\begin{definition}\label{definition:genNakajima}
We say that $G$ is a \define{generalised Nakajima group}
if there exists an ordered basis $\mathcal{B}\subset V^*$ and a sequence of
positive integers $1\leq i_0< i_1< \cdots <i_r\leq n$ $(r< n)$, such that

\begin{enumerate}

\item
$gx_i-x_i\in \Bbbk\langle x_{1},\ldots, x_{i-1}\rangle$ for all
$g\in G$ and for all $i=1,\ldots,n$.

\item
$G$ is generated by the set $P_1\cup P_2\cup\cdots \cup P_r$ as a group
where, for $k=1,\ldots,r$, the subgroup $P_k$ of $G$ is generated by
$$
\{g\in G\mid \text{if}\; gx_i\neq x_i \;\text{then}\; i_{k-1}< i\leq
i_{k} \;\text{and}\; \beta_g\leq i_{k-1}\}.
$$
\end{enumerate}

We call $\mathcal{B}$ a \define{generalised Nakajima basis}
and $(i_0, i_1,\ldots,i_r)$ a \define{generalised Nakajima
sequence} with respect to the generalised Nakajima basis $\mathcal{B}$.
\end{definition}

Note that if $G\subseteq \GL_{\Bbbk}(V)$ is a Nakajima group
and if
$\mathcal{B}=\{x_1,\ldots,x_n\}\subset V^*$ is a Nakajima basis \cite[Definition 8.0.4]{CampbellWehlauModularInvThy11}  for the action of $G$
on $V$, then $G$ is generalised Nakajima group with respect to the basis
$\mathcal{B}$ and the generalised Nakajima sequence $(i_0,\ldots, i_{n-1})=(1,\ldots,n)$.
Hence Nakajima groups are generalised Nakajima groups.
An example of Stong, studied
in~\cite[Example~4.5]{ShankWehlauTransferModular1999}
and~\cite[Example~8.1]{CampbellWehlauModularInvThy11}
is an example of a generalised Nakajima group
that is not a Nakajima group; see Example~\ref{example:stong} for details.

We prove the following theorem:

\begin{theorem}
\label{theorem:genNakajimaHilbertIdealCI}
Let $G$ be a generalised Nakajima group with respect to a basis $\mathcal{B}=\{x_1,\ldots,x_n\}$ of $V^*$. Then
the Hilbert ideal $\hilbertIdeal_{G,S}$ is a complete intersection. Moreover, it can be generated
by a set of homogeneous generators $\{f_1,\ldots,f_n\}$ such that $f_j\in \Bbbk[x_1,\ldots,x_j]$; $\deg(f_j)\leq |G|$ and $x_j^{\deg(f_j)}$ is a
term of $f_j$ for each $j$.
\end{theorem}

\begin{corollary}
Let $G$ be a generalised Nakajima group.
If $S^G$ is a direct summand of $S$ as an $S^G$-module, then $S^G$ is a
polynomial ring.
\end{corollary}

The corollary establishes a conjecture of  R.~J. Shank and D.~L.  Wehlau,
reformulated by A.~Broer, which we now describe.
The \define{trace homomorphism} $\Trace^G : S\longrightarrow S^G$
given by $f \mapsto \sum_{g\in G}gf$ is a homomorphism of
$S^G$-modules.
Shank and Wehlau~\cite[Conjecture~1.1]{ShankWehlauTransferModular1999}
conjectured that $S^G$ is a polynomial ring if and only if
$\image \Trace^G$ is a principal ideal of $S^G$.
Broer~\cite[Corollary~4]{BroerDirectSummandProperty2005} showed that
$\image \Trace^G$ is a principal ideal of $S^G$
if and only if
$S^G$ is a direct summand of $S$ as an $S^G$-module.
Putting these together, the conjecture can be stated as follows:
$S^G$ is a polynomial ring if and only if
$S^G$ is a direct summand of $S$ as an $S^G$-module; what is unknown is the
`if' direction.
Note that, if $S^G$ is a direct summand of
$S$ as an $S^G$-module, and $\hilbertIdeal_{G,S }$ is a complete
intersection, then $S^G$ is a polynomial ring;
see~\cite[Proposition~2.1]{BroerInvThyAbelianTransvGps2010}
for details.
Hence to prove the conjecture for a class of representations $G \to
\GL(V)$, it suffices to show that $\hilbertIdeal_{G,S }$ is a complete
intersection.

In~\cite[Corollary~1.2]{BroerInvThyAbelianTransvGps2010},
Broer proved the above conjecture for abelian $p$-groups,
by showing that if $G$ is an abelian $p$-group generated by
non-diagonalizable pseudo-reflections on $V$
then $\hilbertIdeal_{G,S}$ is a complete
intersection~\cite[Theorem~3.1]{BroerInvThyAbelianTransvGps2010}.
J.~Elmer and M.~Sezer~\cite{ElmerSezerLocallyFiniteDeriv2016}
have obtained results in this direction for some classes of representations
of $p$-groups, using locally finite derivations.
Motivated by this, we try to find classes of representations for which
the Hilbert ideal is a complete intersection.

The bound on the degree of the generators was conjectured by
H.~Derksen and
G.~Kemper~\cite[Conjecture~3.2.6(b)]{DerksenKemperComputationalInvariantThy2ed2015}
(Conjecture~3.8.6(b) in the 2002 edition of the book) for arbitrary finite
groups $G$.
If $|G|$ is invertible in $\Bbbk$,
$\hilbertIdeal_{G,S}$ is generated by
homogeneous elements of degree at most $|G|$;
see~\cite[Lemma~3.2.1]{DerksenKemperComputationalInvariantThy2ed2015} for
a proof by D.~Benson, and references to other proofs. This conjecture is
known to hold in some modular situations also; see the paragraph
after~\cite[Conjecture~3.2.6(b)]{DerksenKemperComputationalInvariantThy2ed2015}.

To prove Theorem~\ref{theorem:genNakajimaHilbertIdealCI},
we introduce a property $P(G',W)$
in Definition~\ref{definition:propertyPGW},
where $G'$ is a subgroup of $G$ and $W$ is a
subspace of $V^{G'}$.
We show inductively
that for a generalised Nakajima group $G$,
there exists a sequence of subgroups
$G_1\subseteq G_2\subseteq \cdots\subseteq G_r=G$
such that for each $k=1,\ldots,r$, the property $P(G_k, W^{(i_k)})$ holds
for some $\Bbbk$-subspace $W^{(i_k)}$ of $V^{G_k}$.

\subsection*{Acknowledgements}
We thank Anurag Singh and Peter Symonds for introducing us to the
Shank-Wehlau conjecture.
The computer algebra systems~\cite{M2}, \cite{Singular}
and~\cite{sagemath}
provided valuable assistance in studying examples.

\section{Hilbert ideal}
\label{section:prelim}

\begin{definition}%
[\protect{\cite[p.~406]{BroerInvThyAbelianTransvGps2010}}]
Let $W \subseteq V^G$ be a subspace. The
\define{Hilbert ideal of $G$ in $S$ with respect to $W$}, denoted
$\hilbertIdeal_{G,S,W}$,
is the $S$-ideal $((W^\perp S) \cap S^G)S$,
where $W^\perp = \ker (V^* \to W^*)$,
i.e., the subspace of linear forms on $V$ that vanish on $W$.
\end{definition}

Note that the Hilbert ideal $\hilbertIdeal_{G,S}$ of $G$ in $S$
is $\hilbertIdeal_{G,S,0}$.
\begin{proposition}%
[\protect{\cite[Lemma~2.2 and its proof]{BroerInvThyAbelianTransvGps2010}}]
\label{proposition:Broer22}
Let $W \subseteq V^G$ be a subspace of codimension $s$.
If $\hilbertIdeal_{G,S,W}$ is a complete intersection generated by
homogeneous invariants $f_1, \ldots, f_s$, then there exist homogeneous
invariants $f_{s+1}, \ldots, f_n$ such that
$\hilbertIdeal_{G,S}$ is the complete intersection ideal generated by
$f_1, \ldots, f_n$.
\end{proposition}

A homogeneous $S$-ideal $I$ is said to be a \define{complete
intersection} if it is minimally generated by $\height(I)$ elements.
A note about terminology:
Broer~\cite[p.~406]{BroerInvThyAbelianTransvGps2010} calls
$\hilbertIdeal_{G,S,W}$ a complete intersection if it is generated by
$\codim_V(W)$ homogeneous invariants. Since $G$ is finite, the map $S^G \to
S$ is finite; further $S^G$ is a normal domain. Hence
$\codim_V(W) = \dim_\Bbbk W^\perp = \height (W^\perp S)
=\height ((W^\perp S) \cap S^G) =\height (\hilbertIdeal_{G,S,W})$.

\begin{notation}
Let $1 \leq i \leq n$.
Let $W^{(i)}$
be the $n-i$ dimensional subspace
$W^{(i)} =\Bbbk\langle v_{i+1},\ldots, v_n\rangle$ and $S^{(i)} =
\Sym(W^{(i)})^\perp=\Bbbk[ x_{1},\ldots, x_{i}]$. 
\end{notation}

\begin{discussionbox}
\label{discussionbox:restriction}
We collect some remarks about the action of $G$ on $V$ and on $V^*$.

\begin{asparaenum}

\item 
$x_1 \in S^G$, since $gx_i -x_i \in \Bbbk \langle x_1, \ldots,
x_{i-1} \rangle$ for every $g \in G$ and $1 \leq i \leq n$.

\item
Let $G'$ be a subgroup of $G$ and let $W$ be a subspace of $V$ with $W\subseteq V^{G'}$.
Then the natural map $V \to V/W \simeq (W^\perp)^*$
is $G'$-equivariant; hence the inclusion map $W^\perp \subseteq V^*$ is
$G'$-equivariant. 
Now additionally suppose that $W = W^{(r ) }$ for some $r$.
With respect to the given basis, suppose that
\[
g = \begin{bmatrix} A & 0 \\ B & 1 \end{bmatrix} \in \GL(V)
\]
where $A$ is an ${r \times r}$ lower-triangular unipotent matrix,
$B$ a ${(n-r) \times r}$ matrix and the $0$ and $1$ denote matrices of
$0$s or $1$ of appropriate size. Then, denoting transposes by $(-)^t$,
\[
(g^{-1})^t =
\begin{bmatrix} (A^{-1})^t & -(BA^{-1})^t\\ 0  & 1 \end{bmatrix}
\in \GL(V^*).
\]
Hence
\[
(g^{-1})^t = (A^{-1})^t \in \GL((W^\perp)^{*})
\qquad \text{and} \qquad
g = A \in \GL(W).
\qedhere
\]

\item
It follows from the above discussion that
$W^{(\beta_g)}\subseteq \{v\in V\mid gv=v\}$, for every $g\in G$.

\end{asparaenum}
\end{discussionbox}

\begin{lemma}
\label{lemma:completeIntersection}
Let $G'$ be a subgroup of $G$ such that $W^{(j_0)}\subseteq V^{G'}$ for
some $j_0\in\{0,\ldots, n-1\}$.
Then for all $j\geq j_0$,
$$
\hilbertIdeal_{G', S, W^{(j)}} =
\hilbertIdeal_{G', S, W^{(j_0)}} + (F_{j_0+1},\ldots, F_j)
$$
where
$F_i \in
((W^{(i)})^\perp S)\cap S^{G'})\setminus ((W^{(i-1)})^\perp S\cap S^{G'})$ for
$i=j_0+1,\ldots,j$.

\end{lemma}
\begin{proof}
 Note that $W^{(j)}\subseteq V^{G'}$ for all $j\geq j_0$.
For all $j> j_0$, the ideal
$(W^{(j)})^\perp S\cap S^{G'}/(W^{(j-1)})^\perp S\cap S^{G'}$ is of height
one in the polynomial ring $S^{G'}/(W^{(j-1)})^\perp S\cap S^{G'}$~\cite[Lemma~2.2 and its proof]{BroerInvThyAbelianTransvGps2010}.
Hence we have
$$
(W^{(j)})^\perp S\cap S^{G'}=(W^{(j-1)})^\perp S\cap S^{G'}+(F_j)
$$
for some $F_j\in (W^{(j)})^\perp S\cap S^{G'}$.
Hence we have $\hilbertIdeal_{G', S, W^{(j)}}=\hilbertIdeal_{G', S,
W^{(j_0)}}+(F_{j_0+1},\ldots, F_j)$ with 
$F_i\in ((W^{(i)})^\perp S\cap S^{G'})\setminus 
((W^{(i-1)})^\perp S\cap S^{G'})$
for $i=j_0+1,\ldots,j$.
\end{proof}

We recall $\Bbbk$-linear locally finite iterative higher derivations
and some of its properties from \cite{ElmerSezerLocallyFiniteDeriv2016}.
Let $G'$ be a subgroup of $G$ generated by $g\in G'$ such that $\beta_g\leq
\beta$ for some fixed $\beta\in \mathbb{N}$.
Let $\beta < j \leq n$.
Define $\Bbbk$-linear \define{locally finite iterative higher derivations}
$\Delta_{j}^{(l)} : S\longrightarrow S$ as follows.
Let $t$ be a variable and extend the $G'$-action on $S$ to $S[t]$ with $G'$
acting trivially on $t$.
Let $\phi_j$ be the ring map $S\longrightarrow S$ given by
\[
x_i \mapsto
\begin{cases}
x_i & \text{if}\; i\neq j,\\
x_j+t & \text{if}\; i=j.
\end{cases}
\]
Note that $\phi_j$ is $G'$-equivariant for each $j$.
Define $\Delta_j^{(l)}$ by $\phi_j(f) = \sum \Delta^{(l)}_j(f)t^l$.

 We have, from the $G'$-equivariance of $\phi_j$,
 $$ \sum g(\Delta^{(l)}_j(f))t^l= g(\phi_j(f)) = \phi_j(gf)  = \sum (\Delta^{(l)}_j(gf))t^l$$
for every $g\in G'$. Hence if $f\in S^{G'}$, then $\Delta^{(l)}_j(f)\in S^{G'}$
for all $l\in \mathbb{N}$ and for all $\beta< j\leq n$.
\begin{lemma}\label{lemma:InsideSmallRing}Let $G'$ be a subgroup of $G$ and $W^{(j_0)}\subseteq V^{G'}$ for some $j_0\in\{0,\ldots, n-1\}$. Then for all $j\geq j_0$,
$$\hilbertIdeal_{G', S, W^{(j)}}=\hilbertIdeal_{G', S, W^{(j-1)}}+(f_j)$$
such that $f_j\in S^{(j)}$ is homogeneous, $\deg(f_j)\leq |G'|$ and
$x_j^{\deg(f_j)}$ is a term of $f_j$.
\end{lemma}

\begin{proof}
We have $\hilbertIdeal_{G', S, W^{(j)}}=\hilbertIdeal_{G', S,
W^{(j-1)}}+(F_j)$ with $F_j\in (W^{(j)})^\perp S\cap S^{G'}\setminus
(W^{(j-1)})^\perp S\cap S^{G'}$, by Lemma~\ref{lemma:completeIntersection}.
Since the action of $G$ is graded on the polynomial ring $S$, we may assume
that $F_j$ is homogeneous.
By the choice of $F_j$ it has the least degree among all homogeneous
elements in $(W^{(j)})^\perp S\cap S^{G'}\setminus (W^{(j-1)})^\perp S\cap
S^{G'}$.
Write $F_j=\sum_{\alpha}f_{j, \alpha}z^{\alpha}$ where $f_{j,\alpha}\in
S^{(j)}$, $\alpha=(\alpha_{j+1},\ldots, \alpha_n)\in \naturals^{n-j}$ and
$z^{\alpha}$ denotes the monomial $x_{j+1}^{\alpha_{j+1}}\cdots
x_n^{\alpha_n}$.
Let $\Delta^{(\alpha)}:S\longrightarrow S$ denotes the differential
operator
$\Delta^{(\alpha_{j+1})}_{j+1}\circ\cdots\circ\Delta^{(\alpha_n)}_{n}$.
Write $S(F_j) = \{\alpha\in \naturals^{n-j} \mid f_{j,\alpha}\neq 0\}$.
Label the elements of $S(F_j)$ as $\alpha^{(i)}, 0 \leq i \leq m$
so that $\alpha^{(0)}>\alpha^{(1)}>\cdots>\alpha^{(m)}$, where
$>$ is the reverse lexicographic order on $\naturals^{n-j}$
with
$(0, 0, \ldots, 0,  1 ) >
(0, 0, \ldots, 1,  0 ) >  \cdots >
(0, 1, \ldots, 0,  0 ) >
(1, 0, \ldots, 0,  0 )$.

Define $F_{jk}$ inductively by setting $F_{j0}=F_j$ and
$$
F_{jk}=F_{j,k-1}-\Delta^{(\alpha^{(k-1)})}(F_{j,k-1})z^{\alpha^{(k-1)}}.
$$
Note that $F_{jm}\in S^{(j)}$ and $F_j-F_{jm}\in \hilbertIdeal_{G', S,
W^{(j-1)}}$ since $\Delta^{(\alpha^{(k)})}(F_{j,k})\in (W^{(j-1)})^\perp
S\cap S^{G'}$ by minimality of $\deg(F_j)$ and induction on $k$.
Set $f_j=F_{jm}$.
Note that $f_j=f_{j, 0}$.
Since the $G'$-orbit product of $x_j$ is an
element of $(W^{(j)})^\perp S\cap S^{G'}\setminus (W^{(j-1)})^\perp S\cap
S^{G'}$, by definition of $F_j$ it must have
$x_j^{\deg(F_j)}$ as a term.
Now, by the minimality of $\deg F_j$,
we see that $\deg(f_j)=\deg(F_j)\leq |G'|$.
\end{proof}

\begin{remark}
From Lemma~\ref{lemma:InsideSmallRing}, it is clear that if $G'$ is a
subgroup of $G$ and $W^{(j_0)}\subseteq V^{G'}$ for some $j_0\in\{0,\ldots,
n-1\}$, then for all $j\geq j_0$, the Hilbert ideal $\hilbertIdeal_{G', S,
W^{(j)}}$ of $G'$ with respect to $W^{(j)}$,  is an extended ideal in $S$
from the ring inclusion $S^{(j)}\hookrightarrow S$.
\end{remark}

\section{Generalised Nakajima groups}
We shall continue with the notations from previous section.

\begin{lemma}
\label{lemma:equalIdeal}
Let $G'$ be a subgroup of $G$ such that $x_{j+1},\ldots,x_n\in S^{G'}$.
Then $\hilbertIdeal_{G', S, W^{(j)}} = \hilbertIdeal_{G', S^{(j)}}S$.
\end{lemma}

\begin{proof}
Since
$\hilbertIdeal_{G', S, W^{(j)}} \supseteq \hilbertIdeal_{G', S^{(j)}}S$, it
suffices to prove the other inclusion.
Suppose that $\hilbertIdeal_{G', S, W^{(j)}}=(F_1, \ldots, F_{s})$ where
$F_i\in (W^{(j)})^\perp\cap S^{G'}$.
Write $F_i = \sum_{{\alpha}} f_{i, {\alpha}}
z^{{\alpha}}$ where
${\alpha} = (\alpha_{j+1},\ldots, \alpha_n)\in \naturals^{n-j}$,
$f_{i,{\alpha}}\in S^{(j)}$
and $z^{{\alpha}}$ denotes the monomial
$x_{j+1}^{\alpha_{j+1}}\cdots x_n^{\alpha_n}$.
By definition of $G'$, $g(z^{{\alpha}}) =
z^{{\alpha}}$ for all $g\in G'$ and for all
monomials $z^{{\alpha}}$.
Since
$F_i\in (W^{(j)})^\perp S\cap S^{G'}$, it follows that
$f_{i,{\alpha}}\in (S^{G'})_+$.
Hence $\hilbertIdeal_{G', S, W^{(j)}}$ can be generated by
$\{f_{i,{\alpha}}\mid i=1,\ldots,s\}\subseteq (S^{(j)})^{G'}_+$.
Therefore $\hilbertIdeal_{G', S, W^{(j)}} \subseteq
\hilbertIdeal_{G', S^{(j)}}S$.
\end{proof}

\begin{definition}\label{definition:propertyPGW}
For any subgroup $G'\leq G$ and for any subspace $W\subseteq V^{G'}$, we
say that the \define{property $P(G', W)$ holds}
if
\begin{enumerate}
\item
 $\hilbertIdeal_{G', S, W}$ is a complete intersection ideal generated by
$\{F_1,\ldots , F_j\} \subseteq W^\perp S\cap S^{G'}$, where
$j = \codim_V(W)$.

\item
 $\hilbertIdeal_{G', S, W}=\hilbertIdeal_{G', \Sym(W^\perp)}S$.
\end{enumerate}
\end{definition}

\begin{remark}
Suppose $G$ is a generalised Nakajima group, $\mathcal{B}\subseteq V^*$ a generalised Nakajima basis
and $(i_0,\ldots, i_r)$ a generalised Nakajima sequence  with respect to $\mathcal{B}$. Let $\mathcal{B}^*=\{v_1,\ldots,v_n\}\subseteq V$ be the dual basis of $\mathcal{B}$. Note that from Definition~\ref{definition:genNakajima}  and Discussion~\ref{discussionbox:restriction}, $G$ can be described in the following equivalent way:

\begin{enumerate}
\item $gv_i-v_i\in \Bbbk\langle v_{i+1},\ldots, v_n\rangle \text{ for all }g\in G, \text{ and for all }i=1,\ldots,n$.
\item $G$ is generated by the set $P_1\cup\cdots\cup P_r$ where the subgroup $P_k $ for $k=1,\ldots, r$ is generated by $g\in G$ such that $W^{(i_{k-1})}\in V^{P_k}$ and $g(v_{i})-v_i\in
W^{(i_{k-1})}\setminus W^{(i_{k})}$  for all $1\leq i \leq i_{k-1}$.
\end{enumerate}
\end{remark}

\begin{lemma}
\label{lemma:SectCI}
Let $G$ be a generalised Nakajima group and
$\mathcal{B}=\{x_1,\ldots,x_n\}$ be a generalised Nakajima basis of $V^*$,
and $(i_0, i_1,\ldots,i_r)$ a generalised Nakajima sequence with respect to
$\calB$.
Let $1 \leq k \leq r$.
Let $G_k$ be the subgroup of $G$ generated by the set
$P_1\cup P_2\cup\cdots \cup P_k$.
Then $P(G_k, W^{(j)})$ holds for all $j\geq i_k$.
\end{lemma}

\begin{proof}
Note that $W^{(i_0)}\subseteq V^{G_1}$ and
$\hilbertIdeal_{G_{1}, S, W^{(i_0)}} =
((x_1,\ldots,x_{i_0})S\cap S^{G_1})S \subseteq
(x_1,\ldots,x_{i_0})S\subseteq \hilbertIdeal_{G_1, S^{(i_0)}}S$;
hence
$\hilbertIdeal_{G_{1}, S, W^{(i_0)}} =
\hilbertIdeal_{G_1, S^{(i_0)}}S =
(x_1,\ldots,x_{i_0})S$.
We thus see that $P(G_1, W^{(i_0)})$ holds.
By Lemmas~\ref{lemma:completeIntersection} and~\ref{lemma:equalIdeal},
$P(G_1, W^{(j)})$ holds for all $j=i_1,\ldots,n$.
Suppose that $P(G_k, W^{(j)})$ holds for all $j=i_k,\ldots,n$.
We shall prove that $P(G_{k+1}, W^{(j)})$ holds for all
$j=i_{k+1},\ldots,n$.
Note that $G_{k+1}$ and $G_k$ act identically on $S^{(j)}$ for
all $j\leq i_k$.
Hence
\begin{equation}
\label{equation:HilbertIdealonSmallPolynomialRing}
\hilbertIdeal_{G_{k+1}, S^{(i_k)}}S =
\hilbertIdeal_{G_{k}, S^{(i_k)}}S
\stackrel{\star}{=}
\hilbertIdeal_{G_k, S, W^{(i_k)}} \supseteq
\hilbertIdeal_{G_{k+1}, S, W^{(i_k)}} \supseteq
\hilbertIdeal_{G_{k+1}, S^{(i_k)}}S,
\end{equation}
where the equation $(\star)$ follows from $P(G_k, W^{(i_k)})$.
Therefore, $\hilbertIdeal_{G_{k+1}, S, W^{(i_k)}} =
(F_1,\ldots,F_{i_k})$
where
$\{F_1,\ldots,F_{i_k}\} \subseteq (W^{(i_{k})})^\perp S\cap S^{G_{k}}$.
Now by Lemmas~\ref{lemma:completeIntersection} and~\ref{lemma:equalIdeal},
$P(G_{k+1}, W^{(j)})$ holds for all $j=i_{k+1},\ldots,n$.
\end{proof}
\begin{proof}[Proof of%
Theorem~\protect{\ref{theorem:genNakajimaHilbertIdealCI}}]
Let $G$ be a generalised Nakajima group with a generalised Nakajima basis
$\mathcal{B}$ and let $(i_0,\ldots, i_r)$ be a generalised Nakajima sequence with
respect to $\mathcal{B}$. By Lemma~\ref{lemma:SectCI}, $P(G, W^{(j)})$ holds for all $j=i_r, \ldots, n$.
Hence $\hilbertIdeal_{G,S}$ is a complete intersection ideal.

Since $\hilbertIdeal_{G_{1}, S, W^{(i_0)}}=\hilbertIdeal_{G_1, S^{(i_0)}}S=(x_1,\ldots,x_{i_0})S$, applying
Lemma~\ref{lemma:InsideSmallRing} with $G'=G_1$ and using Lemma~\ref{lemma:equalIdeal},
we get  $\hilbertIdeal_{G_1, S^{(i_1)}}S=(f_1,\ldots, f_{i_1})$ with
$f_i\in S^{(i)}$ homogeneous, $\deg(f_i)\leq |G_1|$ and $x_i^{\deg(f_i)}$
a term of $f_i$ for each $i=1,\ldots, i_1$.

By~\eqref{equation:HilbertIdealonSmallPolynomialRing}
 we have $\hilbertIdeal_{G_{k+1}, S^{(i_k)}}S =
\hilbertIdeal_{G_{k}, S^{(i_k)}}S=\hilbertIdeal_{G_k, S, W^{(i_k)}}=
\hilbertIdeal_{G_{k+1}, S, W^{(i_k)}}$.
By induction on $k$, we can assume $\hilbertIdeal_{G_{k}, S^{(i_k)}}S=\hilbertIdeal_{G_{k+1}, S, W^{(i_k)}}$
 is generated by a homogeneous set of generators $\{f_1,\ldots,f_{i_k}\}$ such that $f_i\in S^{(i)}$, $\deg(f_i)\leq |G_k|$
 and $x_i^{\deg(f_i)}$ is a term of $f_i$ for each $i=1,\ldots,i_k$.
Now the induction step follows from
 Lemma~\ref{lemma:InsideSmallRing} and Lemma~\ref{lemma:equalIdeal}
by taking $G'=G_{k+1}$.
\end{proof}

\section{Examples}

\begin{example}[Stong]
\label{example:stong}
See~\cite[Example~4.5]{ShankWehlauTransferModular1999}
or~\cite[Example~8.1]{CampbellWehlauModularInvThy11}.
Let $\Bbbk = \bbF_p^3$, with $\bbF_p$-basis $\{1, \omega, \nu \}$.
Let $\rho, \sigma, \tau \in \GL_3(\Bbbk )$, whose action on
$S = \Bbbk[x_1, x_2, x_3]$ is given by
\[
\begin{array}{lll}
\begin{array}{l}
\rho(x_1) = x_1 \\
\rho(x_2) = x_2+x_1 \\
\rho(x_3) = x_3
\end{array}
&
\begin{array}{l}
\sigma(x_1) = x_1 \\
\sigma(x_2) = x_2 \\
\sigma(x_3) = x_3+x_1
\end{array}
&
\begin{array}{l}
\tau(x_1) = x_1 \\
\tau(x_2) = x_2 + \omega x_1\\
\tau(x_3) = x_3 + \nu x_1.
\end{array}
\end{array}
\]
Let $G$ be the subgroup of $\GL_3(\Bbbk)$ 
generated by $\rho$, $\sigma$ and $\tau$.
Then $G$ is a generalised Nakajima group with
generalised Nakajima sequence $(i_0, i_1 ) = (1,3 )$.
It is not a Nakajima group.
\end{example}

\begin{example}
Let $G = \ints/p\ints$; let $\sigma$ be a generator of $G$.
Let $V_2 = \Bbbk \langle v_1, v_2 \rangle$ be a two-dimensional 
vector-space with a $G$-action given by $\sigma(v_1) = v_1-v_2, \sigma(v_2)
= v_2$. Let $m \geq 2$ be an integer. By $mV_2$, we mean the 
$V_2^{\oplus m}$ considered as a representation of 
$G \times \cdots \times G$ ($m$-fold product).
Note that $(V_2^{\oplus m})^* = \Bbbk\langle x_1, y_1, \ldots, x_m,
y_m\rangle$.
Embed $G$ diagonally in $G \times \cdots \times G$.
Then the action of $G$ on $S := \Sym (mV_2)^*$
is given by $\sigma(x_i) = x_i$ and $\sigma(y_i) = y_i+x_i$ for each $i$.
By ordering the variables $x_1, \ldots, x_m, y_1,
\ldots, y_m$, and taking $(i_0,i_1) = (m,2m)$ we see that $G$ is a
generalised Nakajima group.
The description of $S^G$ is given
in~\cite[Theorem~7.4.1 and Remark~7.4.2]{CampbellWehlauModularInvThy11}.
From this we see that $S^G$ is not a polynomial ring and that
$\hilbertIdeal_{G,S} = (x_1, \ldots, x_m, y_1^p, \ldots, y_m^p)$, which is
a complete intersection. 
\end{example}

\begin{example}
Let $\Bbbk=\bbF_3$ and $V=\bbF_3^4$. Let $\sigma \in \GL(V^*)$
be the element that fixes $x_1$ and $x_2$
and maps
$x_3\mapsto x_3+x_2$
and
$x_4\mapsto x_4+x_2+x_1$.
Let $\tau$ be the element that
fixes $x_1$ and $x_4$
and maps
$x_2\mapsto x_2+x_1$
and
$x_3\mapsto x_3+x_1$.
Let $G\subseteq \GL(V^*)$ be generated by $\sigma$ and $\tau$.
Using~\cite{M2} and~\cite{Singular}, we check
that $G$ has order $27$,
that it is not generated by pseudo-reflections,
and
that $S^G$ is not a polynomial ring (which must be the case since $G$ is
not generated by its pseudo-reflections).
However, the Hilbert ideal $\hilbertIdeal_{G,S}=(x_1,x_2^3, x_3^9, x_4^3-x_3^3)$ is a complete intersection
ideal.
\end{example}

\ifreadkumminibib
\bibliographystyle{alphabbr}
\bibliography{kummini}
\else

\def\cfudot#1{\ifmmode\setbox7\hbox{$\accent"5E#1$}\else
  \setbox7\hbox{\accent"5E#1}\penalty 10000\relax\fi\raise 1\ht7
  \hbox{\raise.1ex\hbox to 1\wd7{\hss.\hss}}\penalty 10000 \hskip-1\wd7\penalty
  10000\box7}

\fi %

\end{document}